\documentclass[10pt]{article}

\usepackage[alphabetic]{amsrefs}
\usepackage{amsmath,amssymb,enumerate}
\usepackage[all,cmtip,arrow]{xy}
\usepackage[applemac]{inputenc}

\def \a{{\mathfrak a}}
\def \al{\alpha}

\def \CA{{\cal A}}
\def \CB{{\cal B}}
\def \CC{{\cal C}}

\def \CF{{\cal F}}
\def \CG{{\cal G}}
\def \CH{{\cal H}}
\def \CI{{\cal I}}
\def \CJ{{\cal J}}
\def \CO{{\cal O}}

\def \CR{{\cal R}}

\def \CZ{{\cal Z}}
\def \coim{{\rm coim}\,}
\def \coker{{\rm coker\,}}
\def \df{\ \begin{array}{c} _{\rm def}\\ ^{\displaystyle =}\end{array}\ }
\def \doubleto{\to\hspace{-8pt}\to}
\def \ds{\displaystyle}
\def \e{\emph}
\def \eps{\varepsilon}

\def \F{{\mathbb F}}
\def \For{{\rm For}}

\def \Ga{\Gamma}
\def \Hom{{\rm Hom}}
\def \Id{{\rm Id}}
\def \im{{\rm im}\,}

\def \Map{{\rm Map}}
\def \mathqed{\tag*{$\square$}}
\def \Mod{{\rm Mod}}
\def \ol{\overline}
\def \opp{{\rm opp}}
\def \ph{\varphi}
\def \Q{{\mathbb Q}}

\def \df{\ \begin{array}{c} _{\rm def}\\ ^{\displaystyle =}\end{array}\ }
\def \prf{{\bf Proof: }}
\def \qed{\ifmmode\eqno \square
	\else\noproof\vskip 12pt plus 3pt minus 9pt \fi}
 	\def\noproof{{\unskip\nobreak\hfill\penalty50\hskip2em\hbox{}%
     \nobreak\hfill $\square$\parfillskip=0pt%
     \finalhyphendemerits=0\par}}
\def \Set{{\rm Set}}
\def \setminus{\begin{picture}(18,10)\put(4,6)
                {\line(2,-1){10}}\end{picture}}
\def \Spec{{\rm Spec}\,}

\def \Z{{\mathbb Z}}           
\def \({\left(}
\def \){\right)}
\def \={{\ =\ }}
     
\newcommand{\tto}[1]{\stackrel{#1}{\longrightarrow}}

\newtheorem{theorem}{Theorem}[subsection]

\newtheorem{lemma}[theorem]{Lemma}
\newtheorem{corollary}[theorem]{Corollary}
\newtheorem{proposition}[theorem]{Proposition}

\newtheorem{defi}[theorem]{Definition}
\newenvironment{definition}[0]{\begin{defi}\rm}
{\end{defi}}

\begin{document}

\pagestyle{myheadings} \markright{BELIAN CATEGORIES}

\title{Belian categories}
\author{Anton Deitmar\\ \ \\ 
Far East J. Math. Sci. 70, Issue 1, 1-46 (2012)}

\date{}
\maketitle

{\bf Abstract.} We instal homological algebra, including derived functors, on certain non-additive categories like categories of pointed CW-complexes, modules of monoids or sheaves thereof.
We apply this theory to monoid schemes and sheaves on them, compare the result with the base change and prove several structural theorems.

{\bf MSC: 18A05}, 11G25, 14A15, 14F05, 18F20, 55N30.


{\scriptsize
\tableofcontents}

\section*{Introduction}

In this paper we extend homological algebra to non-additive categories.
Such endeavor has been taken up before, see
\cites{Christensen, DP, EilMoore,  Ina, Keller, Schneiders}, however, non of these fit the situation envisaged by the current paper, which is motivated by the activity around $\F_1$-geometry or non-additive geometry \cites{Manin, KOW, Soule,F1,Haran,Connes,ConnesMonHyp,ConnesNotGeom,Lorscheid}.
The ``common core'' of all these theories seems to be the theory of monoidal schemes as developed in  \cites{Kato,F1}.
In this paper we develop homological algebra over categories which behave like abelian categories, but lack additivity.
The model cases we have in mind are categories of modules over monoidal schemes.
The first part of the paper is devoted to the foundational work on homological algebra for non-additive categories which includes other examples as well, like categories of pointed topological spaces or sheaves thereof.

In the second part of the paper we verify the conditions in the context of sheaves over monoid schemes.
The first part is more algebraic in nature, the second part is more geometric.
We prove some of the results one might expect, like vanishing of cohomology in degrees above the dimension or that cohomology can be computed using flabby resolutions.
Finally, the quite useful compatibility with base change is proved. This allows one to compute the $\Z$-lift of cohomology by means of ordinary Zariski-sheaf cohomology.

I thank Alexander Schmidt for useful remarks on the contents of this paper.

\section{Belian categories}
\subsection{Definition}
\begin{definition}
A category is called \e{balanced} if every morphism which is a monomorphism as well as an epimorphism, already has an inverse, i.e., is an isomorphism.
For example, the category of groups is balanced, but the category of rings is not, since the map $\Z\to\Q$ is an epi- and a monomorphism, but not an isomorphism.
\end{definition}

\begin{definition}
A category $\CC$ is \e{pointed} if it has an object $0$ such that for every object $X$ the sets $\Hom(X,0)$ and $\Hom(0,X)$ have exactly one element each.
The zero object is uniquely determined up to unique isomorphism.
In every set $\Hom(X,Y)$ there exists a unique morphism which factorizes over the zero object, this is called the zero morphism.
In a pointed category it makes sense to speak of kernels and cokernels.
Kernels are always mono and cokernels are always epimorphisms.
\end{definition}

Assume that kernels and cokernels always exist. Then every kernel is the kernel of its cokernel and every cokernel is the cokernel of its kernel.
For a morphism $f$ let $\im(f)=\ker(\coker(f))$ and $\coim(f)=\coker(\ker(f))$.
If $\CC$ has enough projectives, then the canonical map $\coim(f)\to\im(f)$ has zero kernel and if $\CC$ has enough injectives, then this map has zero cokernel.

\begin{definition}
A \e{belian category} is a balanced pointed category $\CB$ which   \begin{itemize}
\item contains finite products, kernels and cokernels, and
\item has the property that
every morphism with zero cokernel is an epimorphism.
\end{itemize}
\end{definition}

Every abelian category is belian.

The second axiom says that a morphism with zero cokernel is an epimorphism and consequently a monomorphism with zero cokernel is an isomorphism.
However, not every morphism with zero kernel is a monomorphism.
Also, not every epimorphism is a cokernel.

Every monomorphism is a kernel.

\begin{definition}
A morphism is called \e{strong}, if the natural map from $\coim(f)$ to $\im(f)$ is an isomorphism.
A belian category 
\end{definition}

Note that this map has zero cokernel, therefore is an epimorphism, so $f$ being strong is equivalent to the map $\coim(f)\to\im(f)$ being injective. 

Kernels and cokernels are strong. Monomorphisms are strong.
If $A\stackrel{f}{\to}B\stackrel{g}{\to}C$ is given with $g$ being strong and $gf=0$, then the induced map $\coker(f)\to C$ is strong. Likewise, if $f$ is strong and $gf=0$, then the induced map $A\to\ker g$ is strong.
A map is strong if and only if it can be written as a cokernel followed by a kernel.

A belian category $\CB$ is called \e{strong}, if every morphism in $\CB$ is strong.
If $\CB$ is any belian category, the subcategory $\CB^{\rm strong}$ that has the same objects, but ony the strong morphisms of $\CB$ is again a belian category but this time a strong one. 

Note that in a belian category, although one cannot add morphisms, one can ``add'' morphisms from direct sums thanks to the universal property of direct sums: Suppose given two morphisms $\ph_i\colon M_i\to N$, $i=1,2$. Then there exists a unique morphism
$$
\ph_1\oplus\ph_2\colon M_1\oplus M_2\ \to\ N
$$
such that 
$
\xymatrix{
M_i\ar[r]& M_1\oplus M_2\ar[r]^{\ \ \ \ \ph_1\oplus \ph_2}& N
}
$
 equals $\ph_i$ for $i=1,2$.

{\bf Example.} 
The simplest example of a belian category is the category $\Set_0$ of pointed sets.
Objects are pairs $(X,x_0)$ where $X$ is a set and $x_0\in X$ is an element. A morphism $\ph\in\Hom((X,x_0),(Y,y_0))$ is a map $\ph\colon X\to Y$ with $\ph(x_0)=y_0$. Any singleton $(\{ x_0\},x_0)$ is a zero object. The kernel of a morphism $\ph\colon X\to Y$  is the inverse image $\ph^{-1}(\{ y_0\})$ of the special point and the cokernel is $Y/\ph(X)$, where the image $\ph(X)$ is collapsed to a point.
The product is the cartesian product and the coproduct is the disjoint union with the special points identified.
A morphism $\ph\in\Hom((X,x_0),(Y,y_0))$ is strong if and only if $\ph$ is injective outside $\ph^{-1}(\{ y_0\})$.

Other examples include the category of pointed simplicial sets, pointed $CW$-complexes, or the categories of sheaves of such.

If $\CB$ is a belian category, then for $X,Y\in\CB$ the set $\Hom_\CB(X,Y)$ is a pointed set, the special point being the zero morphism.

\subsection{Complexes}
\begin{definition}
In a belian category a sequence of morphisms, 
$$
\dots\tto{\ }M^i\tto{d^i}M^{i+1}\tto{d^{i+1}}\dots
$$
is called a \e{complex} if $d^{i+1}\circ d^i=0$ for every $i$.
In that case there is an induced morphism $\im d^i\to\ker d^{i+1}$ which is a monomorphism since the maps $\im(d^i)\to M^{i+1}$ and $\ker(d^{i+1})\to M^{i+1}$ are monomorphisms.
We call the complex \e{exact}, if this morphism is an isomorphism.
For a given complex let
$$
H^i(M^\bullet)\df \coker\(\im d^i\to\ker d^{i+1}\)\ \in\ \CB
$$
be the \e{cohomology} of the complex $M^\bullet$.
Then the cohomology is zero if and only if the complex is exact.
A complex is called a \e{strong complex} if every differential $d^i$ is strong.
\end{definition}

Let $\CB$ be a belian category and let $\CC(\CB)$ be the category of complexes over $\CB$.
Morphisms in $\CC(\CB)$ are morphisms $f\colon X\to Y$ of complexes, i.e., $f$ is a sequence $f^i\colon X^i\to Y^i$ of morphisms is $\CB$ such that every square
$$
\xymatrix{
{X^i}\ar[r]\ar[d]^{f^i}
	&{X^{i+1}}\ar[d]^{f^{i+1}}\\
{Y^i}\ar[r]	
	&{Y^{i+1}}
}
$$
is commutative.

Let $\CC_+(\CB)$ be the full subcategory of complexes $Y$ which are bounded below, i.e., $Y^i=0$ for $i <<0$.
Further $\CC_-(\CB)$ denotes the subcategory of complexes which are bounded above and finally let $\CC_b(\CB)=\CC_+(\CB)\cap \CC_-(\CB)$ be the category of bounded complexes.

\subsection{Pull-backs and push-outs}
\begin{definition}
Let $\CB$ be a category.
An object $I\in \CB$ is called \e{injective} if for every monomorphism $M\hookrightarrow N$ the induced map $\Hom(N,I)\to \Hom(M,I)$ is surjective.
Dually, an object $P\in\CB$ is called \e{projective} if it is injective in the opposite category $\CB^{\opp}$ with reversed arrows.
We say that $\CB$ \e{has enough injectives} if for every $A\in\CB$ there exists a monomorphism $A\hookrightarrow I$, where $I$ is an injective object.
Likewise, we say that $\CB$ \e{has enough projectives} if for every $A\in\CB$ there is an epimorphism $P\doubleto A$ with $P$ projective, or, equivalently, if $\CB^\opp$ has enough injectives.
\end{definition}

\begin{lemma}\label{epi}
Let $\CB$ be a category and let
$$
\xymatrix{
A\ar[r]^{f'}\ar[d]
	& X \ar[d]\\
{B}\ar[r]_{f}
	&{Y}
}
$$
be a cartesian square in $\CB$.
\begin{itemize}
\item If $f$ is a monomorphism, then so is $f'$.
\item If $\CB$ contains enough projectives and $f$ is an epimorphism,
then $f'$ is an epimorphism.
\item If $\CB$ is belian and contains enough injectives and enough projectives, and if $f$ is strong, then $f'$ is strong.
\end{itemize}

Likewise, let 
$$
\xymatrix{
{A}\ar[r]^{h}\ar[d]
	&{B}\ar[d]\\
{C}\ar[r]^{h'}
	&{P}
}
$$
be co-cartesian.
\begin{itemize}
\item If $h$ is an epimorphism, then so is $h'$.
\item If $\CB$ contains enough injectives and $h$ is a monomorphism, then $h'$ is a monomorphism.
\item If $\CB$ is belian and contains enough injectives, $h$ is a strong morphism then $h'$ is strong and there is an isomorphism $C/\ker h \cong C/\ker h'$.
\end{itemize}
Finally, if $\CB$ is belian and contains enough injectives, then every monomorphism is a kernel.
In particular, a morphism is strong if and only if it can be written as a cokernel followed by an injection.
\end{lemma}

Note that an epimorphism is in general not a cokernel.

\prf
Assume the first situation and let $\al,\beta$ be two morphisms $Z\to A$ with $f'\al= f'\beta$.
We have to show $\al=\beta$.
Since $fg'\al =fg'\beta$ and $f$ is injective, we have $g'\al= g'\beta$.
The square being cartesian implies $\al=\beta$ as claimed.
For the second assertion,
let $\al\colon P\to X$ be an epimorphism with $P$ projective.
The resulting morphism $P\to Y$ can be lifted to $B$, giving a commutative square
$$
\xymatrix{
{P}\ar[r]\ar[d]
	&{X}\ar[d]\\
{B}\ar[r]
	&{Y.}
}
$$
Since the original square was cartesian, the epimorphism $P\to X$ factorizes as $P\to A\stackrel{f'}{\to} X$, hence $f'$ is an epimorphism.
We postpone the proof of the third property till later.

The first two assertions for co-cartesian squares follow by reversing the arrows.
Before proving the third, we first prove the final remark that every monomorphism is a kernel.
So assume the belian category $\CB$ to contain enough injectives.
Let $f: A\hookrightarrow B$ be a monomorphism in $\CB$.
Let $B/A$ denote the cokernel of $f$ and let $K$ be the kernel of $B\to B/A$.
We have the following diagram,
$$
\xymatrix{
{A}\ar@{^{(}->}[r]^{f}\ar@{^{(}->}[d]
	&{B}\ar@{->>}[r]
		&{B/A}\\
{K}\ar@{^{(}->}[ur]\ar@{->>}[d]\\
{C,}
}
$$
where $C$ is the cokernel of the natural map $A\to K$, which necessarily must be injective.
We have to show that $C$ is zero. Then the second axiom implies that the map from $A$ to $K$ is surjective as well and as the category is balanced, it is an isomorphism.
We want to show that the diagram
$$
\xymatrix{
{K}\ar@{^{(}->}[r]\ar@{->>}[d]
	&{B}\ar@{->>}[d]\\
{C}\ar[r]^{0}
	&{B/A}
}
$$
is co-cartesian.
Once this is shown, the claim follows, as by the above, the zero morphism in the bottom then is injective, hence $C=0$.
So assume given arrows $C\to Z$ and $B\to Z$ which become the same on $K$. Consider the diagram
$$
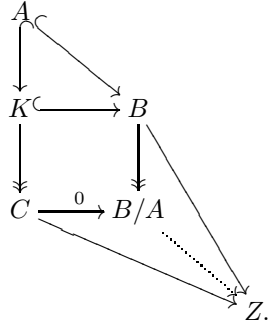
\xymatrix{
{A}\ar@{^{(}->}[d]\ar@{^{(}->}[dr]\\
{K}\ar@{^{(}->}[r]\ar@{->>}[d]
	&{B}\ar@{->>}[d]\ar[ddr]\\
{C}\ar[r]^{0}\ar[drr]
	&{B/A}\ar@{.>}[dr]\\
&&{Z.}
}
$$
The solid arrow diagram commutes.
As $A\to K\to C$ is zero, then $A\hookrightarrow B\to Z$ is zero, so a unique dotted arrow exists, making the triangle $B, B/A, Z$ commutative.
As the rest of the diagram commutes, this implies that the triangle $C,B/A,Z$ also commutes, i.e., the entire diagram is commutative, which implies that the square indeed is co-cartesian and it follows that $f$ is a kernel indeed.

Now for the third assertion on co-cartesian diagrams. 
Let $K$ be the kernel of $h$ and write $A/K$ for its cokernel.
Likewise let $K'$ be the kernel of $h'$ and $C/K'$ its cokernel.
We get the solid arrow diagram
$$
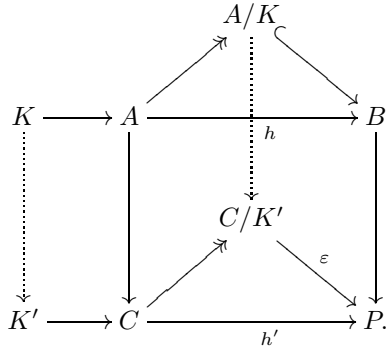
\xymatrix{
&&{A/K}\ar@{^(->}[dr]\ar@{.>}[dd]\\
{K}\ar[r]\ar@{.>}[dd]
	&{A}\ar@{->>}[ur]\ar[rr]_{\ \ \ \ h}\ar[dd]
	&&{B}\ar[dd]\\
&&{C/K'}\ar[dr]^\eps\\
{K'}\ar[r]
	&{C}\ar@{->>}[ur]\ar[rr]_{\ \ \ \ h'}
		&&{P.}
}
$$
The dotted arrows are implied by the kernel and cokernel properties.
It is easy to see that $P$ also is the cofiber product of $C/K'$ and $B$ over $A/K$.
Therefore, $\eps$ is injective, and so $h'$ is strong.
Replacing $C/K'$ with $C/K$ the same arguments work, so we get $C/K'\equiv C/K$.

We finally have the ingredients to also prove the third claim on catesian squares.
Let $K$ and $K'$ be kernels to $f$ and $f'$ and let $B/K$ as well as $A/K'$ denote their cokernels.
We get the solid arrow diagram
$$
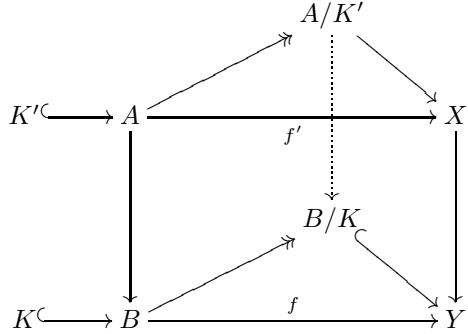
\xymatrix{
&&&{A/K'}\ar[dr]\ar@{.>}[dd]\\
{K'}\ar@{^(->}[r]
	&{A}\ar[rrr]_{f'}\ar@{->>}[rru]\ar[dd]
	&&&{X}\ar[dd]\\
&&&{B/K}\ar@{^(->}[dr]\\
{K}\ar@{^(->}[r]
	&{B}\ar[rrr]^{f}\ar@{->>}[rru]
	&&&{Y}
}
$$
There exists a unique dotted arrow making the diagram commutative by the cokernel property of $A/K'$.
We claim that the right hand square is also cartesian.
For this let $Z$ be an object with arrows to $X$ and $B/K$ giving the same arrow to $Y$.
As $\CB$ has enough projectives, we can replace $Z$ with a projective object $P$.
We get the diagram
$$
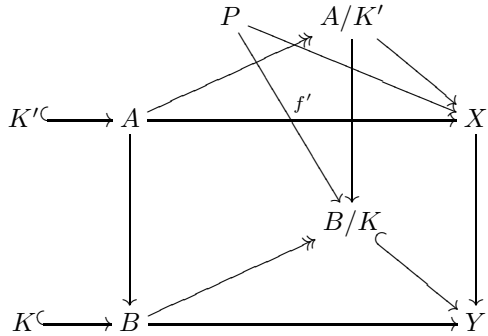
\xymatrix{
&&{P}\ar[rrd]\ar[ddr]
	&{A/K'}\ar[dr]\ar[dd]\\
{K'}\ar@{^(->}[r]
	&{A}\ar[rrr]^{f'}\ar@{->>}[rru]\ar[dd]
	&&&{X}\ar[dd]\\
&&&{B/K}\ar@{^(->}[dr]\\
{K}\ar@{^(->}[r]
	&{B}\ar[rrr]\ar@{->>}[rru]
	&&&{Y}
}
$$
As $P$ is projective, there exists an arrow from $P$ to $B$ making the diagram commutative.
Next the cartesian property of the original square implies the existence of an arrow from $P$ to $A$ and the resulting arrow from $P$ to $A/K'$ shows that the right hand square is indeed cartesian as claimed.
Then the injectivity of the arrow from $B/K'$ to $Y$ implies the injectivity of the arrow from $A/K'$ to $X$, the latter therefore is a kernel and the map $f'$ is strong.
\qed

\begin{definition}
Let $\CB$ be a category which contains fiber-products and has enough projectives.
Let $Y$ be an object in $\CB$.
On the class of morphisms $h\colon X\to Y$ with target $Y$ we define an equivalence relation as follows.
We say that $(h,X)\sim (h',X')$ if there exists a commutative diagram
$$
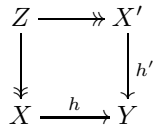
\xymatrix{
{Z}\ar@{->>}[d]\ar@{->>}[r]
	&{X'}\ar[d]^{h'}\\
{X}\ar[r]^{h}
	&{Y}
}
$$
where the arrows emanating at $Z$ are epimorphisms.
One has to check that this indeed is an equivalence. The only problem is transitivity.
For this assume $(h,X)\sim (h',X')$ and $(h',X')\sim (h'',X'')$.
This means that we have the solid arrows in the following diagram,
$$
\xymatrix{
{Z''}\ar@{.>}[r]\ar@{.>}[d]
	&{Z'}\ar@{->>}[r]\ar@{->>}[d]
		&{X''}\ar[d]^{h''}\\
{Z}\ar@{->>}[r]\ar@{->>}[d]
	&{X'}\ar[r]^{h'}\ar[d]^{h'}
		&{Y}\\
{X}\ar[r]^{h}
	&{Y.}
}
$$
Let $Z''$ be the fiber-product so that the upper left square is cartesian.
Then by the last lemma the dotted arrows are epimorphisms and so are the arrows $Z''\to X$ and $Z''\to X''$.
This proves that $\sim$ is an equivalence relation.
Equivalence classes are called \e{generalized elements}.
They are a useful technical tool as they allow proofs by diagram chase as the following lemma shows.
By abuse of notation, we will write $y\in Y$ for a class $y=[h]$ with target object $Y$.
If $\al:Y\to Z$ is a morphism, we write $\al(y)$ for the class of the morphism $\al h$.
We will write $|Y|$ for the class of generalized elements of the object $Y$.
Note that if the category is pointed, then $|Y|$ is a pointed class.
\end{definition}

\begin{lemma}
Assume that $\CB$ is a balanced pointed category with fibre products and  enough injectives and projectives.
Then a sequence
$$
{A}\tto{\al}{B}\tto{\beta}{C}
$$
is exact if and only if the induced sequence
$$
{|A|}\tto{\al}{|B|}\tto{\beta}{|C|}
$$
of generalized elements is exact.

\end{lemma}

\prf Assume the sequence is exact.
Then $\beta\al=0$.
Let $b\in |B|$ with $\beta(b)=0$.
Let $K$ be the kernel of $\beta$.
Then $\al$ factorizes over $K\cong\im(\al)$, and as $\CB$ has enough injectives, the map $A\to\im(\al)=K$ is surjective.
Since $\beta(b)=0$, $b$ factorizes over $K$ as well.
We have the following diagram of solid arrrows,
$$
\xymatrix{
{P}\ar@{->>}[rr]\ar@{.>}[dd]
	&&{X}\ar[dl]\ar[dd]^{b}\ar[ddr]^{0}\\
&{K}\ar@{^(->}[dr]	\\
{A}\ar@{->>}[ur]\ar[rr]
	&&{B}\ar[r]
	&	{C,}
}
$$
where $P$ is any projective cover of $X$.
By projectivity, the map $P\to K$ lifts to $A$, giving the dotted arrow which is the searched for pre-image of $b$.

For the converse direction assume $\beta\al=0$ and the condition on elements.
We have the diagram,
$$
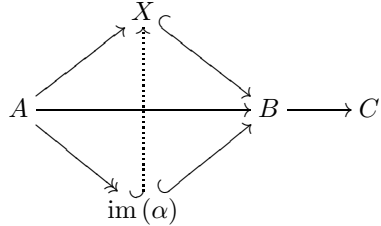
\xymatrix{
&{X}\ar@{^(->}[dr]\\
{A}\ar[ur]\ar[rr]\ar@{->>}[dr]
	&&{B}\ar[r]
		&{C}\\
&{\im(\al)}\ar@{^(.>}[uu]\ar@{^(->}[ur]
}
$$
The condition applied to the element $[K\hookrightarrow B]$ yields a map $P\to A$ and a surjection $P\twoheadrightarrow K$ making the diagram
$$
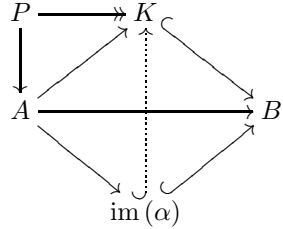
\xymatrix{
{P}\ar@{->>}[r]\ar[d]
	&{K}\ar@{^(->}[dr]\\
{A}\ar[ur]\ar[rr]\ar[dr]
	&&{B}\\
&{\im(\al)}\ar@{^(->}[ur]\ar@{^(.>}[uu]
}
$$
commutative.
As the map $P\twoheadrightarrow K$ is onto, so is $\im(\al)\to K$, which therefore is an isomorphism as $\CB$ is balanced.
\qed

\subsection{Snake Lemma}
\begin{definition}
Let $\CB$ be a belian category with fibre products and let $X\in\CB$ be an object.
A generalized element $x\in|X|$ is called an \e{atom} if for every strong morphism $f:X\to Y$ with $f(x)\ne 0$ and every $x'\in|X|$ one has
$$
f(x)=f(x')\quad\Rightarrow\quad x=x'.
$$
An \e{atomic class} is a class $J$ of atoms which is stable under strong morphisms, such that with $|X|_J=|X|\cap J$ we have that a sequence $A\to B\to C$ is exact if and only if the induced sequence $|A|_J\to|B|_J\to|C|_J$ is exact.
If an atomic class exists, we will consider it fixed and call any atom in $J$ an \e{admissible atom}.

We say that a belian category $\CB$ \e{admits an atomic class}, if it is closed under fibre products and there exists an atomic class in $\CB$.
\end{definition}

\begin{definition}
A morphism $\ph:X\to Y$ in a belian category is called \e{pseudo-isomorphism} if it has trivial kernel and cokernel. 
A pseudo-isomorphism is onto, but not necessarily injective.
In that case we also say that $X$ is a \e{pseudo-isomorphic cover} of $Y$.

Note that in this situation if one of $X$, $Y$ is zero, then so is the other.
\end{definition}

\begin{lemma}(Snake Lemma)
Let $\CB$ be a category which is belian, has enough injectives and projectives and admits an atomic class.
Given a  commutative diagram with exact rows
$$
\xymatrix{
&{X_1}\ar[r]^{g_1}\ar[d]^{f_1}
		&{X_2}\ar[r]^{g_2}\ar[d]^{f_2}
			&{X_3}\ar[r]\ar[d]^{f_3}
				 &0\\
{0}\ar[r]
	&{Y_1}\ar[r]^{h_1}
		&{Y_2}\ar[r]^{h_2}
			&{Y_3,}
}
$$
called a \e{snake diagram},
the induced sequences
$$
\ker(f_1)\to\ker(f_2)\to\ker(f_3)
$$
and
$$
\coker(f_1)\to\coker(f_2)\to\coker(f_3)
$$
are exact.
If $h_1$, $f_2$ and $g_2$ are strong, then there is a natural strong morphism
$\delta\colon \ker(f_3)\to\coker(f_1)$ such that the whole sequence is exact.
\end{lemma}

\begin{corollary}\label{Cor1.4.4}
Given a snake diagram with $h_1$ and $f_2$ strong but $g_2$ not strong, one can replace $g_2$ with the cokernel of $g_1$ to obtain a snake diagram  where $X_3\tto {f_3}Y_3$ is replaced with $\coker(g_1)=\tilde X_3\tto{\tilde f_3}Y_3$ and there exists a pseudo-isomorphism $\ker(\tilde f_3)\to \ker(f_3)$. For the modified diagram there exists a snake morphism $\delta$.
\end{corollary}

\begin{corollary}[Weak snake]\label{Cor1.4.5}
Given a snake diagram in which only $h_1$ is strong.
After the replacement of Corollary \ref{Cor1.4.4} one obtains a the following diagram with exact rows
$$
\xymatrix{
&&\ker(f_1)\ar[r]
&\ker(f_2)\ar[r]^\al
&\ker(\delta)\ar@{=}[d]
\\
\coker(f_1)
&\coker(f_2)\ar[l]
&\coker(f_3)\ar[l]
&\ker(\tilde f_3)\ar[l]_\delta
&\ker(\delta).\ar@{_(->}[l]
}
$$
The map $\delta$ is not necessarily strong.
\end{corollary}

\prf
The first corollary is clear. The second will be proved after the proof of the Lemma.
The exactness of the two sequences is obtained by a standard verification.
To construct $\delta$,
extend the diagram as follows:
$$
\xymatrix{
&{X_1}\ar[r]^{t}\ar[d]^{\Id}
		&{Z}\ar@{->>}[r]^{s}\ar@{^(->}[d]^{l}
			&{\ker(f_3)}\ar[r]\ar@{^(->}[d]^{k}
				 &0\\
&{X_1}\ar[r]^{g_1}\ar[d]^{f_1}
	&{X_2}\ar[r]^{g_2}\ar[d]^{f_2}
		&{X_3}\ar[r]\ar[d]^{f_3}
			&0\\
{0}\ar[r]
	&{Y_1}\ar[r]^{h_1}\ar@{->>}[d]^{k'}
		&{Y_2}\ar[r]^{h_2}\ar@{->>}[d]^{l'}
			&{Y_3}\ar[d]^{\Id}\\
&{\coker(f_1)}\ar[r]^{s'}
	&{Z'}\ar[r]^{t'}
		&{Y_3.}
}
$$
Here $Z$ is the kernel of $f_3g_2$.
As $k$ is a kernel and $g_2$ a cokernel, $Z$ happens to be the fibre product of $\ker(f_3)$ and $X_2$ over $X_3$.
Next $Z'$ is the cokernel of $h_1f_1$.
As $h_1$ is a kernel and $k'$ a cokernel, $Z'$ happens to be the cofibre product of $\coker(f_1)$ and $Y_2$ over $Y_1$.
By Lemma \ref{epi}, $s$ is an epimorphism.
Likewise, $s'$ is a monomorphism.
The morphism $t$ is the fibre product of $g_1$ and the zero map from $X_1$ to $\ker(f_3)$.
We claim that the first row is exact.
Since $st=0$ it remains to show that $t$ is surjective on $\ker(s)$.
Now $g_1$ is surjective on $\ker(g_2)$.
Replacing $g_1$ by $\ker(g_2)$ amounts to the same as assuming that $g_1$ is injective. It suffices to prove the claim under that assumption.
Indeed, then $t$ is the kernel of $s$.
To see this, let $W\stackrel{w}{\to}Z$ be a morphism with $sw=0$.
We shall show that $w$ factorizes uniquely over $t$.
The induced arrow $W\to X_3$ is zero, therefore there is a unique morphism $r\colon W\to X_1$ such that the solid arrow diagram
$$
\xymatrix{
{W}\ar[r]^w\ar[d]_r
	&{Z}\ar[d]^l\\
{X_1}\ar[r]\ar@{.>}[ur]^t
	&{X_2}
}
$$
is commutative.
We have to show that it remains commutative when $t$ is inserted.
We have two morphisms $w,tr\colon W\to Z$ with $lw=ltr$ and $ sw=str=0$.
By the universal property of the fibre product $Z$ it follows that $w=tr$, hence the diagram commutes.
Further, by construction the morphisms $k$ and $k'$ are strong, so by Lemma
\ref{epi} the morphisms $l$ and $l'$ are strong.

In the last row the morphism $t'$ is the cofibre product of $h_2$ and zero.
The exactness of this row
follows from the previous part by reversing all arrows.
So the rows are exact.
Consider the morphism $\eps = l'f_2 l\colon Z\to Z'$. It satisfies $\eps t=s'k'f_1 =0$ and $t'\eps = f_3ks =0$.
Since $s$ is the cokernel of $t$, there exists a unique map from $\ker(f_3)$ to $Z'$ making the ensuing diagram commutative.
Next $s'$ is the kernel of $t'$, so the map from $\ker(f_3)$ to $Z'$ factors uniquely over $\coker(f_1)$ giving the desired map  $\delta\colon \ker(f_3)\to\coker(f_1)$ such that $\eps=s'\delta s$.
We claim that $\delta$ is strong if $f_2$ is.
In this case the map $\eps$ is strong.
Consider the commutative diagram with exact rows
$$
\xymatrix{
&{X}\ar[r]^{t}\ar[dr]^{0}
	&{Z}\ar[r]^s\ar[dr]^0\ar[d]^\eps
		&{\ker f_3}\ar[r]&0\\
{0}\ar[r]&{\coker f_1}\ar[r]^{s'}
	&{Z'}\ar[r]^{t'}
		&{Y_3.}
}
$$
Firstly,the induced morphism $\delta_1 :\ker f_3\to Z'$ such that $\eps=\delta_1 s$ is strong, as there are natural isomorphisms $\coim \delta_1\cong\coim\eps$ and $\im\delta_1\cong\im\eps$ identifying the natural map $\coim\delta_1\to\im\delta_1$ with $\coim\eps\to\im\eps$ which is an isomorphism.
Similarly, the natural map $\delta$ such that $\delta_1=s'\delta$ is strong.
It remains to show that the sequence
$$
\ker(f_2)\tto\al\ker(f_3)\tto\delta\coker(f_1)\tto\beta\coker(f_2)
$$
is exact.
From the construction of $\delta$ it follows that $s'\delta\al=0$ and since $s'$ is mono, we get $\delta\al=0$.

Next let $0\ne x\in|\ker(\delta)|$ be an admissible atom.
We have to show that there exists an admissible atom $u$ in $|\ker(f_2)|$ such that $\al(u)=x$.
Pick an admissible pre-image $0\ne z\in |Z|$ under $s$.
We have to show that $z$ can be chosen such that $l(z)\in|\ker(f_2)|$.
We have $l'f_2l(z)=0$, which means $f_2(l(z))\in\ker l'$ and the latter equals the image of $\ker(k')$ by Lemma \ref{epi}.
By exactness we have $\ker(k')=\im(f_1)$ so that there exists an admissible atom $w\in |X_1|$
with $f_2(l(z))=f_2(g_1(w))$.
By admissiblity, $l(z)$ is an atom, as $f_2$ is strong, it follows either $f_2(l(z))=0$, or $l(z)=g_1(w)$.
As the first case is what we want, we deal with the second now.
If $l(z)=g_1(w)$, then by commutativity of the diagram, $t(w)=z$, but then $x=s(z)=s(t(w))=0$, which we have excluded.
Together we have shown exactness of the sequence at $\ker(f_3)$.

We now show $\beta\delta=0$.
As $\beta$ is induced by $h_1$ it suffices to show $s'\delta=0$, but as $s'\delta s=0$ and $s$ is onto, this is clear.

Next let $x\in|\ker(\beta)|$ be an admissible atom, so $x\in |\coker(f_1)|$.
Pick an admissible pre-image $y$ in $|Y_1|$.
Then $l'h_1(y)=0$. As the map $Y_2\to\coker(f_2)$ factors over $l'$ it follows that $h_1(y)$ maps to zero in $\coker(f_2)$, hence there is an admissible $v\in |X_2|$ with $f_2(v)=h_1(y)$.
This element maps to zero in $Y_3$, so it comes from $\ker(f_3)$ and we have found a pre-image under $\delta$.
This finishes the proof of the snake lemma.
\qed

{\bf Proof of Corollary \ref{Cor1.4.5}.}
One checks that the only points of the proof, where the strongness of $f_2$ is used, is that $\delta$ is strong and $\ker(\delta)/\im(\al)=0$.
\qed

As an application we will show the existence of a long exact cohomology sequence attached to a short exact sequence of complexes.

\begin{lemma}
We assume that $\CB$ has enough injectives and projectives and admits an atomic class.
Let
$$
S\equiv\quad 0\to E\stackrel{e}\to F\stackrel{f}{\to}G\to 0
$$
be an exact sequence of complexes over the belian category $\CB$, where the map $e$ and the complex $F$ are strong.
Then there is a pseudo-isomorphic cover $\tilde H^p(S)\doubleto H^p(G)$ and a commutative diagram,
$$
\xymatrix{
H^p(F)\ar[r]^{\tilde H^p(f)}\ar[dr]_{H^p(f)} & \tilde H^p(S)\ar@{->>}[d]\\
& H^p(G),
}
$$
together with an exact sequence
$$
\cdots\to H^p(F)\tto{\al} \tilde H^p(S)\stackrel{\delta}{\to}H^{p+1}(E)\to H^{p+1}(F)\to\cdots.
$$
Then the connection maps $\delta$ are strong.
 The connection map is functorial in the sense that for a commutative diagram of complexes,
 $$
\xymatrix{
S\equiv\ar[d]
&0\ar[r] 	
&E\ar[r]^e\ar[d]	
&F\ar[r]\ar[d]	
&G\ar[r]\ar[d]			
&0\\
T\equiv
&0\ar[r]
&X\ar[r]^\xi
&Y\ar[r]
&Z\ar[r]
&0,
}
$$
where $e$, $\xi$, $F$ and $Y$ are strong, one gets a commutative diagram
$$
\xymatrix{
H^p(F)\ar[r]\ar[d]
&\tilde H^p(S)\ar[r]^\delta\ar[d]
&H^{p+1}(E)\ar[d]\\
H^p(Y)\ar[r]
&\tilde H^p(T)\ar[r]^\delta\
&H^{p+1}(X)
}
$$
for every $p\in\Z$.
The cover  $\tilde H^p$ is a functor from the category of short exact sequences  $S$ to $\CB'$.

If the complex $F$ is not strong, one obtains the following natural  diagram with exact rows,
$$
\xymatrix{
\ker(\delta)\ar@{=}[d]
&H^p(F)\ar[l]_\al
&H^p(E)\ar[l]
\\
\ker(\delta)\ar@{^(->}[r]
&\tilde H^p(S)\ar[r]^\delta
&H^{p+1}(E)\ar[r]
&H^{p+1}(F)\ar[r]
&\tilde H^{p+1}(S)
}
$$
\end{lemma}

\prf
At each stage $p\in\Z$ one gets a  commutative and exact diagram,
$$
\xymatrix{
&{\coker d_E^{p-1}}\ar[d]\ar[r]^{\xi'}
	&{\coker d_F^{p-1}}\ar[d]\ar[r]
		&{\coker d_G^{p-1}}\ar[d]\ar[r]
			&{0}\\
{0}\ar[r]
	&{\ker d_E^{p+1}}\ar[r]^{\xi''}
		&{\ker d_F^{p+1}}\ar[r]
			&{\ker d_G^{p+1},}
}
$$
where $\xi'$ and $\xi''$ are strong.
Replacing $\coker d_G^{p-1}$ with the pseudo-isomorphic cover $\coker(\xi')$ is functorial in the sequence.
The snake lemma gives the desired long exact sequence.
\qed

\subsection{Ascent functors}
\begin{definition}
A functor between belian categories is called \e{strong-exact} if it maps strong exact sequences to exact sequences.
Note that being strong-exact is a weaker condition than being exact.
It is, however, the natural condition in the non-additive context.
\end{definition}

\begin{definition}
Let $\CB$ be a belian category.
An \e{ascent functor}  $A:\CB\to\CC$ is a functor from $\CB$ to an abelian category $\CC$, such that
\begin{itemize}
\item $A$ is faithful and strong-exact, and
\item $A$ maps epimorphisms to epimorphisms.
\end{itemize}
\end{definition}

Note that an ascent functor will preserve the canonical factorization of a strong morphism into a cokernel followed by a kernel, and it will preserve the classes of epi- and monomorphisms.

{\bf Example.} The standard example will be the category $\CB$ of pointed sets and the ascent functor $A:\CB\to\Mod(\Z)$, mapping a pointed set $(X,x_0)$ to the $\Z$-module $\Z[X]/\Z x_0$.
The properties are easily verified.

\begin{lemma}\label{ascentProp}
Let $A$ be an ascent functor on the belian category $\CB$.
For any morphism $f$ in $\CB$ we have canonical maps
\begin{enumerate}[\rm (a)]
\item $A(\coker(f))\cong\coker(A(f))$,
\item $A(\im(f))\cong\im(A(f))$, and
\item $A(\ker(f))\hookrightarrow \ker(A(f))$ which is an isomorphism, if $f$ is strong.
\end{enumerate}
Consequently, for a complex $M^\bullet$ we have a canonical injection
$$
A(H^i(M^\bullet))\ \hookrightarrow\ H^i(A(M^\bullet)),
$$
which is an isomorphism if the complex is strong.
\end{lemma}

\prf
To prove these points, let $f:X\to Y$.
The sequence
$$
0\to\im(f)\to Y\to \coker(f)\to 0
$$
is strong and exact, so it will remain exact after applying $A$, which shows that (a) and (b) imply each other.
To prove (a), consider the diagram
$$
\xymatrix{
&Z\\
A(X)\ar[ur]^0\ar[r]\ar@{->>}[d]
& A(Y)\ar[u]\ar[r]
&A(\coker(f))\ar[r]\ar@{.>}[ul]&0\\
A(\coim(f))\ar@{->>}[r]
&A(\im(f))\ar[u],
}
$$
where all but the dotted arrow are given, the surjections are preserved by $A$. 
Then, as the map from $A(X)$ to $A(\im(f))$ is surjective, it follows that the map from $A(\im(f))$ to $Z$ is zero.
Since $A(\coker(f))$  is the cokernel of $A(\im(f))\hookrightarrow A(Y)$, it follows the unique existence of the dotted arrow making the entire diagram commutative.
This fact then shows that indeed there is a canonical isomorphism $A(\coker(f))\cong\coker(A(f))$.

It remains to prove (c).
For this we first recall that the sequence 
$$
0\to\ker(f)\to X\to\coim(f)\to 0
$$
is strong-exact, therefore remains exact after applying $A$, which implies that $A(\ker(f))$ injects into $A(X)$.
We therefore have all arrows but the dotted one in the diagram
$$
\xymatrix{
A(\ker(f))\ar@{.>}[d]\ar@{^(->}[dr]\\
\ker(A(f))\ar@{^(->}[r]&A(X)\ar[r]^{A(f)}&A(Y)
}
$$
As the map from $A(\ker(f))$ to $A(Y)$ is zero, it follows that a unique dotted arrow exists, as the map from $A(\ker(f))$ to $A(X)$ is injective, the dotted arrow is injective, too.
If $f$ is strong, the sequence $0\to A(\ker(f))\to A(X)\to A(Y)$ is exact.
\qed

\begin{lemma}\label{AscentStrongness}
Let $A$ be an ascent functor on the belian category $\CB$.
Consider a sequence 
$$
\xymatrix{
S\ \equiv&X\ar[r]^g&Y\ar[r]^f&Z
}
$$
 in $\CB$ and assume that the induced sequence
$
{A(S)}
$
is exact.
Then the sequence $S$ is exact in $\CB$.
If $g$ is strong, then $f$ is strong as well.
\end{lemma}

\prf
Assume $A(S)$ is exact. Then $0=A(f)A(g)=A(fg)$ and the faithfulness implies $fg=0$.
So we get a natural map $\ph:\im g\to \ker f$.
Applying $A$, we get a commutative diagram
$$
\xymatrix{
\im A(g)\ar[r]^\cong&A(\im g)\ar[dr]\ar@{^(->}[r]^{A(\ph)}&A(\ker f)\ar@{^(->}[d]\\
&&\ker A(f).
}
$$
As the sequence $A(S)$ is exact, the diagonal arrow is an isomorphism, so then is $A(\ph)$.
As $A$ is faithful, $\ph$ is epi and mono, hence also an isomorphism as $\CB$ is balanced.
So the sequence $S$ is exact.

For the second assertion let $g$ be strong. It suffices to assume that $g$ is the kernel of $f$.
The sequence $X\to Y\to Y/X\to 0$ is strong and exact, therefore the sequence $A(X)\to A(Y)\to A(Y/X)\to 0$ is exact, consequently one has $A(Y/X)\cong A(Y)/A(X)$, and the latter injects into $A(Z)$.
As $A$ is faithful, $Y/X$ also injects into $Z$, so $f$ is strong.
\qed

Note that an ascent functor will in general not preserve products, as the following example shows.
Let $\CB$ be the belian category of pointed sets and let $A:\CB\to \Mod(\Z)$ be the functor that maps a pointed set $(X,x_0)$ to the $\Z$-module $\Z[X]/\Z x_0$.
Then, for two finite sets $X,Y\in\CB$ the $\Z$-module $A(X\times Y)$ is free of dimension $\# X\cdot \# Y -1$, whereas $A(X)\times A(Y)$ has dimension $\#X-1+\# Y -1$.

\subsection{Derived functors}
\begin{definition}
A \e{resolution} of an object $X$ in $\CB$ is a strong exact sequence
$$
0\to X\to I_X^0\to I_X^1\to\cdots .
$$
A functor $F\colon\CB\to\CB'$ between belian categories is called \e{left strong-exact} if  $F$ maps kernels to kernels and for every strong exact sequence
$$
0\to X\stackrel{\ph}{\hookrightarrow} Y\to Z
$$
in $\CB$, the sequence
$$
0\to F(X)\hookrightarrow F(Y)\to F(Z)
$$
is exact in $\CB'$. 
\end{definition}

Let $F\colon \CB\to\CB'$ be left strong-exact and assume that $\CB$ has enough injectives.
Then every ebject $X$ possesses a resolution $0\to X\to I^0\to I^1\to\dots$, where all $I^p$ are injective objects. Such a resolution is called an \e{injective resolution}.
Fix an injective resolution $X\to I_X$ for every $X\in\CB$.

\begin{lemma}
For every morphism $f:X\to Y$ in $\CB$ there exist morphisms $f^0,f^1,\dots$ making the diagram
$$
\xymatrix{
0\ar[r]
& X\ar[r]\ar[d]_f
&I_X^0\ar[r]\ar[d]_{f^0}
&I_X^1\ar[r]\ar[d]_{f^1}&\dots\\
0\ar[r]
& Y\ar[r]
&I_Y^0\ar[r]
&I_Y^1\ar[r] &\dots
}
$$
commutative.
\end{lemma}

\prf
The map $f^0$ exists as $I_Y^0$ is injective.
Next $I_X^0/X$ injects into $I_X^1$, so, by the same reason, $f^1$ exists.
Iteration yields the lemma.
\qed
 
\begin{definition}
For any morphism $f:X\to Y$ in $\CB$ fix a choice of morphisms $f^0,f^1,\dots$ as in the lemma.
For $p=0,1,\dots$ define
$$
R^p_IF(X)\df H^p(F(I_X^\bullet)),
$$
and for every $f:X\to Y$ set
$$
R^p_IF(f)=H^p(f^\bullet),
$$
the induced morphism on cohomology.
\end{definition}

Note that in this general setting, $R^pF$ might not even be a functor.
Only $R^0F$ is, as the next lemma shows.
Later we will show that under additional assumptions $R^pF$ is a functor.

\begin{lemma}
$R^0_IF$ is a functor natural isomorphic with $ F$.
\end{lemma}

\prf
Let $0\to X\to I_X^0\to\cdots$ be the chosen resolution of $X\in \CB$.
Since $F$ is left strong-exact, the sequence $0\to F(X)\hookrightarrow F(I_X^0)\to F(I_X^1)$ is exact.
Therefore there exists a natural functorial isomorphism,
$$
R^0_IF(X)\= H^0F(I_X)\ \cong\ F(X).\qed
$$

\begin{definition}
Let $\CB$ be a belian category.
An \e{injective class} in $\CB$ is a class $\CI$ of injective objects in $\CB$ such that
\begin{itemize}
\item every object of $\CB$ injects into an object in $\CI$, and
\item $\CI$ is closed under finite products.
\end{itemize}
Note that every belian category $\CB$ with enough injectives admits injective classes.
\end{definition}

\begin{definition}
The ascent functor $A$ is said to be \e{$\CI$-injective} if it maps objects in $\CI$ to injective objects.
Here $\CI$ is an injective class.
If we can choose $\CI$ to be the class of all injective objects we simply say that \e{$A$ preserves injectives} and likewise in the projective case.
\end{definition}

{\bf Example.}
Let $\Set_0$ be the category of pointed sets as before.
For a ring $R$ and a pointed set $(M,m_0)$ let $R[M]$ be the free $R$-module generated by $M$ and let $R[M]_0=R[M]/Rm_0$.
Then $A\colon M\mapsto R[M]_0$ from $\Set_0$ to the category of $R$-modules is an ascent functor which is $\CI$-injective for every injective class $\CI$, if $R$ is a field.
Note that this functor indeed is strong-exact but not exact.

\begin{definition}
Let now $F\colon \CB\to\CB'$ be a left strong-exact functor.
An \e{ascent datum} for $F$ is a quadruple $(\CI,A,A', \tilde F)$ 
consisting of an injective class $\CI$ in $\CB$ and an ascent functor $A:
\CB\to\CC$ which is $\CI$-injective, as well as an ascent functor 
$A':\CB'\to\CC'$, and a left-exact functor $\tilde F: \CC\to\CC'$ 
such that on the full subcategory of $\CB$ given by the class $\CI$,  
the functors $A'F$ and $\tilde F A$ from $\CI\subset \CB$ to $\CC'$ are isomorphic.
In other words, the diagram
$$
\xymatrix{
{\CC}\ar[r]^{\tilde F}
	&{\CC'}\\
{\CI}\ar[u]^{A}\ar[r]^{F}
	&{\CB'}\ar[u]_{A'}
}
$$
is commutative up to isomorphy of functors.
\end{definition}

\begin{theorem}\label{Prop1.6.4}
Assume that the left strong-exact functor $F$ is equipped with an ascent datum $(\CI,A,A', \tilde F)$.
Then the $R^pF$ are functors and they depend on the choice of the injective class $\CI$ and the injective resolutions only up to canonical isomorphism.
We have a natural injection
$$
A'R^pF(X)\hookrightarrow R^p\tilde F(A(X)).
$$
On the other hand, if $(\CI,A_1,A_1',\tilde F_1)$ is another ascent datum with the same injective class $\CI$, then this ascent datum will trivially give the same derived functors $R^pF$.
\end{theorem}

\prf
Let first $I_X$ be an arbitrary injective resolution of $X$ and let $A'$ be an ascent functor on $\CB'$.
As
$$
H^i(FI_X)\= \coker(\im(Fd^{i-1})\to\ker(Fd^i)),
$$
Lemma \ref{ascentProp} gives a natural injection
$$
A'(H^i(FI_X))\ \hookrightarrow\ H^i(A'(FI_X)).
$$
We now show that $R_I^pF$ is a functor.
Let $f:X\to Y$ be a morphism. We have to show that 
$R^p_IF(f)$ does not depend on the choice of the maps $f^0,f^1,\dots$.
This proof is completely analogous to the proof that shows independence of resolutions, so we only give the latter.
So suppose there 
are given two injective classes $\CI$ and $\CI_1$ such that $(\CI,A,A', \tilde F)$ and $(\CI_1,A,A', \tilde F)$ form ascent data for $F$.
For $X\in\CB$, choose injective resolutions $I_X$ and  $I_{1,X}$ from $\CI$ and $\CI_1$.
We fix isomorphims of functors
$$
\Phi_I: \tilde FA|_I\tto\cong A'F|_I,
$$
and 
$$
\Phi_{I_1}: \tilde FA|_{I_1}\tto\cong A'F|_{I_1}.
$$

By the injectivity of $I_X$ there is a map extending the identity on $X$,
$$
\xymatrix{
{I_{1,X}}\ar[r]^{\ph}\ar[d]
	&{I_X}\ar[d]\\
{X}\ar[r]^{=}&{X.}
}
$$
Applying the ascent functor $A$ we get a commutative diagram in the abelian category $\CC$,
$$
\xymatrix{
{A(I_{1,X})}\ar[r]^{A(\ph)}\ar[d]
	&{A(I_X)}\ar[d]\\
{A(X)}\ar[r]^{=}&{A(X).}
}
$$
Both columns are injective resolutions.
Therefore the map $A(\ph)$ is unique up to homotopy.
So the map $\tilde F(A(\ph))$ from $\tilde F(A(I_{1,X}))$ to 
$\tilde F(A(I_X))$ is independent of $\ph$ up to homotopy, which implies that the induced map on cohomology,
$$
H^i(\tilde F(A(\ph))): H^i(\tilde F(A(I_{1,X})))\ \to\ H^i(\tilde F(A(I_X)))
$$
is uniquely determined.
The isomorphism $\Phi_I$ induces an isomorphism of complexes $A'FI_X^\bullet\tto\cong \tilde FAI_X^\bullet$ and so an isomorphism of the respective cohomology objects.
The following induced diagram is commutative,
$$
\xymatrix{
{A'(H^iFI_{1,X})}\ar[r]\ar@{^(->}[d]
	&{A'(H^iFI_X)}\ar@{^(->}[d]\\
H^i(A'FI_{1,X})\ar[d]^{\Phi_{I_1}}_\cong\ar[r]&H^i(A'FI_X)\ar[d]^{\Phi_I}_\cong\\
{H^i(\tilde FA(I_{1,X}))}\ar[r]
	&{H^i(\tilde FA( I_X)),}
}
$$
where the horizontal arrows are induced by $\ph$.
As we have learned, the lowest horizontal arrow does indeed not depend on the choice of $\ph$ and by the commutativity of the diagram, so don't the others.
By the faithfulness of ascent we get the uniqueness of $H^i(F(\ph)):H^iFI_{1,X}\to H^iFI_X$.
This then must be an isomorphy by standard arguments.
\qed

\begin{proposition}\label{prop1.6.9}
Let $F\colon\CB\to\CB'$ be a left strong-exact functor on belian categories  equipped with an ascent datum.
Then every injective object of $\CB$ is $F$-acyclic.
\end{proposition}

\prf
Let $J$ be injective and let $J\hookrightarrow I$ be an injection into some $\CI$-injective $I$.
Then, as $J$ is injective, there exists $s:I\doubleto J$ with $si=\Id_J$.
Then $R^p(s)R^pF(i)=R^p(si)=\Id$ and therefore $R^pF(J)$ injects into $R^p(I)$ which is zero for $p>0$, therefore $J$ is acyclic.
\qed

\begin{definition}\label{Def1.7.10}
Let $F:\CB\to\CB'$ be a left strong-exact functor on belian categories equipped with an ascent datum $(A,\CI,A',\tilde F)$.
An \e{acyclic class} for $F$ is a class $\CA$ of objects in $\CB$, such that
\begin{itemize}
\item $\CI\subset\CA$,
\item $A\in\CA$ $\Rightarrow$ $R^pF(A)=0$ for all $p\ge 1$,
\item for every strong exact sequence
$
0\to A\to B\to C\to 0
$ with $A\in\CA$,
the sequence
$$
0\to FA\to FB\to FC\to 0
$$
is exact,
\item if $0\to A\to B\to C\to 0$ is strong and exact, and if $A$ and $B$ belong to $\CA$, then so does $C$.
\end{itemize}
\end{definition}

\begin{theorem}\label{thm1.6.11}
Assume that the left strong-exact functor $F$ is equipped with an ascent datum and an acylic class $\CA$.
Let $0\to X\to A^0\to\dots$ be a resolution of $X$ with $A^p\in\CA$ for all $p\ge 0$.
Then
$$
R^pF(X)\cong H^p(F(A^\bullet)),
$$
i.e., the cohomology can be computed using $\CA$-resolutions.
\end{theorem}

\prf
Choose an $\CI$-resolution
$$
0\to X\to I^0\to I^1\to\dots
$$
such that we get a commutative diagram
$$
\xymatrix{
{0}\ar[r]
	&{X}\ar[r]\ar[d]_{=}
	&	{A^0}\ar[r]\ar@{^(->}[d]
	&		{A^1}\ar[r]\ar@{^(->}[d]
	&			{\cdots}\\
{0}\ar[r]
	&{X}\ar[r]
	&	{I^0}\ar[r]
	&		{I^1}\ar[r]
	&			{\cdots}
}
$$
where the vertical maps can be chosen injective by enlarging $I^j$ is necessary.
Let $(Y^j)$ be the sequence of cokernels so that we get an exact, strong, commutative diagram,
$$
\xymatrix{
&&{0}\ar[d]
	&{0}\ar[d]\\
{0}\ar[r]
	&{X}\ar[r]\ar[d]_{=}
	&	{A^0}\ar[r]\ar[d]
	&		{A^1}\ar[r]\ar[d]
	&			{\cdots}\\
{0}\ar[r]
	&{X}\ar[r]
	&	{I^0}\ar[r]\ar[d]
	&		{I^1}\ar[r]\ar[d]
	&			{\cdots}\\
&{0}\ar[r]
	&{Y^0}\ar[r]\ar[d]
	&	{Y^1}\ar[r]\ar[d]
	&		{\cdots}\\
&&{0}
	&{0}
}
$$
Since $A^p$ and $I^p$ are in $\CA$, so is $Y^p$.
Applying $F$ we obtain a short exact sequence of complexes
$$
0\to F(A)\hookrightarrow F(I)\to F(Y)\to 0.
$$
The corresponding cohomology sequence reads
$$
\tilde H^{i-1}F(Y)\to H^iF(A)\to H^iF(I)\to H^iF(Y).
$$
Both ends are zero, so we get a pseudo-isomorphism in the middle.
However, if $A\hookrightarrow B$ is an injection of complexes such that the induced map $H^p(A)\to H^p(B)$ is a pseudo-isomorphism, then it is an isomorphism, hence
\begin{align*}
H^pF(A)&\cong R^pF(X).\mathqed
\end{align*}

\subsection{Strong derived functors}
\begin{definition}
An object $X$ in $\CB$ is called \e{$F$-acyclic} if $R^iF(X)=0$ for every $i>0$.
\end{definition}

Recall that a functor $F$ on belian categories is called \e{strong}, if it maps strong morphisms to strong morphisms.

\begin{theorem}\label{1.8}
Let $F\colon\CB\to\CB'$ be a left strong-exact functor on belian categories  equipped with an ascent datum.
Assume that $F$ is strong and that $\CB'$ contains enough injectives and projectives and admits an atomic class.
Let $0\to X\to A^0\to A^1\to\cdots$ be a resolution by $F$-acyclics.
Then $R^iF(X)\cong H^i(F(A^\bullet))$, so cohomology can be computed using resolutions by arbitrary acyclics.
\end{theorem}

\prf
We need some lemmas.

\begin{lemma}\label{lem1.7.3}
Under the conditions of the theorem, let 
$$
0\to X\to Y\to Z\to 0
$$
be a strong exact sequence, then for every $p\ge 0$ there exists a pseudo-isomorphic cover $\tilde R^pF(Z)$ of $R^pF(Z)$ and a long exact sequence
$$
0\to F(X)\to F(Y)\to \tilde F(Z)\to R^1F(X)\to R^1F(Y)\to \tilde R^1F(Z)\to\dots
$$
\end{lemma}

\prf
Given a strong exact sequence $0\to X\hookrightarrow Y\to Z\to 0$ in $\CB$ let $I_X$ and $I_Y$ be given $\CI$-resolutions of $X$ and $Z$.
Consider the diagram
$$
\xymatrix{
{}
	&{0}\ar[d]
	&	{0}\ar[d]
	&		{0}\ar[d]\\
{0}\ar[r]
	&{X}\ar@{^(->}[r]\ar[d]
	&	{Y}\ar[r]\ar[d]^{\beta}
	&		{Z}\ar[r]\ar[d]
	&			{0}\\
{0}\ar[r]
	&{I_X^0}\ar[r]^{\al}
	&	{I_X^0\times I_Z^0}\ar[r]^{\gamma}
	&		{I_Z^0}\ar[r]
	&			{0,}
}
$$
where $\al$ is the natural map given by the universal property of the product and the maps $I_X^0\stackrel{\rm id}{\to} I_X^0$ and $I_X^0\stackrel 0\to I_Z^0$.
For the definition of $\beta$ recall that since $I_X^0$ is injective, the map $X\to I_X^0$ extends to $Y\to I_X^0$ and $\beta$ is given by this map and the composition $Y\to Z\to I_Z^0$.
Finally, $\gamma$ is the projection onto the second factor.
The commutativity of the diagram is immediate.
The morphism $\al$ is strong and the diagram is  exact.
Since $\CI$ is an injective class, $I_X^0\times I_Z^0$ lies in $\CI$.

We write $I_Y^0=I_X^0\times I_Z^0$ and extend the diagram by the corresponding cokernels $X', Y', Z'$ to get a commutative  exact diagram such that the left horizontal arrows are strong.
$$
\xymatrix{
{}
	&{0}\ar[d]
	&	{0}\ar[d]
	&		{0}\ar[d]\\
{0}\ar[r]
	&{X}\ar@{^(->}[r]\ar[d]
	&	{Y}\ar[r]\ar[d]
	&		{Z}\ar[r]\ar[d]
	&			{0}\\
{0}\ar[r]
	&{I_X^0}\ar@{^(->}[r]\ar[d]
	&	{I_Y^0}\ar[r]\ar[d]
	&		{I_Z^0}\ar[r]\ar[d]
	&			{0}\\
{0}\ar[r]
	&{X'}\ar@{^(->}[r]\ar[d]
	&	{Y'}\ar[r]\ar[d]
	&		{Z'}\ar[r]\ar[d]
	&			{0}\\
&{0}&{0}&{0.}
}
$$
One uses diagram chase to verify the exactness of this diagram.
We repeat the procedure with the exact sequence $0\to X'\to Y'\to Z'\to 0$.
Iteration leads to a commutative and exact diagram of injective resolutions
$$
\xymatrix{
{}
	&{0}\ar[d]
	&	{0}\ar[d]
	&		{0}\ar[d]\\
{0}\ar[r]
	&{X}\ar@{^(->}[r]\ar[d]
	&	{Y}\ar[r]\ar[d]
	&		{Z}\ar[r]\ar[d]
	&			{0}\\
{0}\ar[r]
	&{I_X}\ar@{^(->}[r]
	&	{I_Y}\ar[r]
	&		{I_Z}\ar[r]
	&			{0.}
}
$$
Applying $F$ to this diagram yields a exact sequence of complexes,
$$
0\to F(I_X)\hookrightarrow F(I_Y)\to F(I_Z)\to 0.
$$
To verify the exactness recall that by construction $I_Y^j$ is the direct product of $I_X^j$ and $I_Z^j$.
For any two objects $A,B$ in $\CB$ the map $A\stackrel{{\rm id}\times 0}\to A\times B\to A$ is an automorphism of $A$. 
Hence the same is true for $F(A)\to F(A\times B)\to F(A)$, so the map $F(A\times B)\to F(A)$ is an epimorphism and $F(A)\to F(A\times B)$ is a monomorphism.
Note that the  sequence $0\to F(I_X)\hookrightarrow F(I_Y)\to F(I_Z)\to 0$ is not strong in general.

To this sequence of complexes we now apply the snake lemma to get a long exact sequence
$$
\cdots \to R^iF(Y)\to \tilde R^iF(Z)\stackrel{\delta}{\to} R^{i+1}F(X)\to R^{i+1}F(Y)\to\cdots
$$
as claimed.
\qed

To finish the proof of the theorem, note that Lemma \ref{lem1.7.3} gives the crucial part of the proof, that the class $\CA$ of all $F$-acyclic objects is an acyclic class in the sense of Definition \ref{Def1.7.10}.
The rest of the properties are clear, so that the present theorem follows from Theorem \ref{thm1.6.11}.
\qed

\begin{proposition}
Let $F\colon\CB\to\CB'$ be a left strong-exact functor on belian categories  equipped with an ascent datum.
If $F$ is strong, and
 $\CB'$ contains enough injectives and projectives and admits an atomic class, then
any injective resolution of $X\in\CB$ computes $R^iF(X)$.
\end{proposition}

\prf
This follows from Proposition \ref{prop1.6.9} together with Theorem \ref{1.8}.
\qed

\section{Pointed modules and sheaves}

\subsection{Definitions}
\begin{definition}
Let $A$ be a commutative monoid. A \e{module} over $A$ is a set $M$ together with an action $A\times M\to M$ sending $(a,m)$ to $am$, satisfying $(ab)m=a(bm)$ and $1m=m$ for all $a,b\in A$ and every $m\in M$.
Let $N\subset M$ be a sub-module, then we define the quotient module $M/N$ by collapsing $N$: as a set, $M/N$ equals $M/\sim$, where $\sim$ is the equivalence relation with the equivalence classes $\{ m\}$, $m\notin N$ and $N$. The module structure is defined by $a[m]=[am]$, where $[m]$ is the class of $m\in M$.

An element $m_0\in M$ is called \e{stationary} if $am=m$ for every $a\in A$.
A \e{pointed module} is a pair $(M,m_0)$ consisting of an $A$-module $M$ and a stationary point $m_0\in M$.
A \e{homomorphism of pointed modules} from $(M,m_0)$ to $(N,n_0)$ is an $A$-module homomorphism $\ph$ with $\ph(m_0)=n_0$.
Let $\Mod_0(A)$ denote the category of pointed modules and their homomorphisms.
The special point $m_0$ of a pointed module $M$ is also denoted by $0_M$ or $0$ if no confusion is likely.
It is called the \e{zero element} of $M$.
\end{definition}

If $M$ is a module over $A$, we define the pointed module $M_0$ to be $M\cup\{ 0\}$, where $0$ is a new stationary point which we choose to be the special point of $M_0$.

The category $\Mod_0(A)$ contains a terminal and initial object, the zero module $\{ 0\}$, also written $0$.
A morphism $\ph:M\to N$ is called zero if $\ph$ factors over zero. This is equivalent to $\ph(M)=\{0_N\}$.

The category $\Mod_0(A)$ contains products and coproducts.
Products are the usual cartesian products and coproducts are given as follows: Let $(M_i)_{i\in I}$ be a family of objects in $\Mod_0(A)$, then the coproduct is
$$
\coprod_{i\in I} M_i\= \left.\bigcup_{i\in I}^\cdot M_i\right/\sim
$$
where the union means the disjoint union of the $M_i$ and the equivalence relation just identifies all zeros $0_{M_i}$ to one.
We also write coproducts as direct sums.

\subsection{Limits}

\begin{proposition}
The category $\Mod_0(A)$ contains direct and inverse limits.
\end{proposition}

\prf
Let $I$ be a small category and $F\colon I\to \Mod_0(A)$ be a functor. Write $M_i$ for $F(i)$, $i\in I$. Define
$$
M\df \coprod_{i\in I}M_i/\sim,
$$
where $\sim$ is the equivalence relation given by $m\sim F(\ph)(m)$ whenever $m\in M_i$ and $\ph\colon i\to j$ is a morphism in $I$.
A straightforward verification shows that $M$ is a direct limit.

Likewise,
$$
N\df \left\{\left.\ x\in\prod_{i\in I}M_i\ \right| \ x_j = F(\ph)(x_i)\ \forall \ph\in\Hom_I(i,j)\ \right\}
$$
is an inverse limit.
\qed

\begin{lemma}
A morphism $\ph\colon X\to Y$ in $\Mod_0(A)$ is an epimorphism if and only if $\ph$ is a surjective map.
\end{lemma}

\prf
Suppose $\ph$ is an epimorphism, then $Y/\im\ph$ is zero, so $\im\ph=Y$, i.e., $\ph$ is surjective.
The rest is clear.
\qed

\subsection{Injectives and projectives, ascent}

\begin{proposition}
The category $\Mod_0(A)$ is a belian category with enough injectives and enough projectives. 
It possesses an atomic class.
\end{proposition}

\prf
It is clear that every morphism with zero cokernel is an epimorphism.
We prove the existence of enough injectives.
For any set $X$ we have an $A$-module structure on the set $\Map(A,X)$ of all maps $\al\colon A\to X$ given by
$$
a\al(b)\= \al(ab).
$$
Further, if $X$ is a pointed set, then $\Map(A,X)$ is a pointed module, the special point being $\al_0$ with $\al_0(a)=x_0$, where $x_0$ is the special point of $X$.
For a given pointed module $M$ we define $I_M$ to be
$$
I_M\df \Map(A,M).
$$
We have a natural embedding $M\hookrightarrow \Map(A,M)$ of $A$-modules given by $m\mapsto \al_m$ with $\al_m(a)=am$.
The theorem will follow if we show that $\Map(A,M)$ is indeed injective.
For this note that for any $A$-module $P$ and any set $X$ there is a functorial isomorphism of $A$-modules
$$
\psi\colon \Map(P,X)\to \Hom_A(P,\Map(A,X))
$$
given by
$$
\psi(\al)(p)(a)\= \al(ap).
$$
The inverse is given by
$$
\psi^{-1}(\beta)(p)\= \beta(p)(1).
$$
Now let $P\hookrightarrow N$ be an injective A-module homomorphism, then for any set $X$ one has the commutative diagram
$$
\xymatrix{
{\Hom_A(N,\Map(A,X))}\ar[r]\ar[d]^{\cong}
&{\Hom_A(P,\Map(A,X))}\ar[d]^{\cong}
\\
{\Map(N,X)}\ar[r]&{\Map(P,X)}
}
$$
The second horizontal map is surjective, therefore the first horizontal map is surjective as well.
For $X=M$ this implies  the first part of the theorem.

For the existence of enough projectives, consider $A$ as a module over itself.
Let $P_M=\bigoplus_{m\in M} A_m^+$ be a direct sum of copies of $A^+$.
Then the pointed module $P_M$ is projective as a straightforward verification shows.
For a given module $M$ define a map
\begin{eqnarray*}
 \ph\colon  P_M &\to & M\\
a\in A_m &\mapsto & am,\\
0 &\mapsto & m_0.
\end{eqnarray*}
Then $\ph\colon P_M\to M$ is the desired surjection.

For an atomic class, we take the class of all generalized elements $[h]$ where $h:M\to X$ and the module $M$ is generated by a single element, i.e.,  there exists $m\in M$ such that $M=Am\cup\{ 0\}$.
It is clear, that this is indeed an atomic class.
\qed

\begin{proposition}\label{2.2.2}
The functor $A$ from $\Mod_0(A)$ to the category of $\Q$-vector spaces,
$$
A(M)\= \Q[M]/Qm_0,
$$
where $m_0$ is the special point, is 
an ascent functor which preserves injectives and projectives.
\end{proposition}

\prf
It is easy to see that $A$ is an ascent functor.
Since every object in the category of $\Q$-vector spaces is injective as well as projective, $A$ preserves these classes of morphisms.
\qed

\subsection{Pointed sheaves}
\begin{definition}
Let $X$ be a monoided space, i.e., a topological space with a sheaf $\CO_X$ of monoids.
A given topological space can be made a monoided space by defining $\CO_X$ to be the constant sheaf $\CO_X(U)=\{ 1\}$.
A \e{pointed sheaf} is a sheaf of pointed $\CO_X$-modules where the restrictions are assumed to preserve the special points.
Let $\Mod_0(X)$ denote the category of pointed sheaves.
\end{definition}

\begin{proposition}
The category   $\Mod_0(X)$ is belian and contains enough injectives.
\end{proposition}

\prf
The zero object is the zero sheaf.
The existence of fiber and cofiber products is a standard sheaf theoretic construction.
To verify the last axiom let $\ph\colon\CF\to\CG$ be a morphism with zero cokernel and let $\CG\ ^{\longrightarrow}_{\longrightarrow}\ \CZ$ be two morphisms such that the induced morphisms from $\CF$ to $\CZ$ agree.
For any $x\in X$ one has the exact sequence of the stalks ${\CF_x}\to\CG_x\to 0$. Therefore $\ph_x$ is an epimorphism and thus the two maps ${\CG_x}\ ^{\longrightarrow}_{\longrightarrow}\ \CZ_x$ agree.
Since this holds for every $x\in X$, the two morphisms $\CG\ ^{\longrightarrow}_{\longrightarrow}\ \CZ$ agree, so $\ph$ is an epimorphism.
\qed

\begin{lemma} The following holds in $\Mod_0(X)$.
\begin{enumerate}[\rm (a)]
\item A morphism $f\colon\CF\to\CG$ is strong if an only if all stalks $f_x\colon\CF_x\to\CG_x$, $x\in X$, are strong.
\item A sequence $\CF\stackrel f\to\CG\stackrel g\to\CH$ is exact if and only if all the sequences at the stalks $\CF_x\stackrel {f_x}\to\CG_x\stackrel {g_x}\to\CH_x$, $x\in X$, are exact.
\end{enumerate}

\end{lemma}

\prf
(a)  A morphism $f$ in a belian category is strong if and only if the induced $\tilde f : \coim f\to\im f$ is an isomorphism.
If $f$ is a morphism in $\Mod_0(X)$, then for every $x\in X$ one has $(\tilde f)_x=\tilde{f_x}$.
Replacing $f$ by $\tilde f$ it therefore suffices to show that $f$ is a monomorphism if and only if all its stalks $f_x$ are.

Let's assume that $f$ is a monomorphism and let $x\in X$.
We have to show that $f_x$ is injective. For this assume $f_x(s_x)=f_x(t_x)$ for some $s_x,t_x\in\CF_x$.
Then there exists an open neighborhood $U$ of $x$ and representatives $s_U,t_U\in\CF_U$ with $f_U(s_U)=f_U(t_U)$ in $\CG(U)$.
We can consider $\CO|_U$ as an $\CO_U$-module, but not a pointed one in general.
To make it pointed we add an extra stationary point $\omega_V$ to $\CO_V$ for every open $V\subset U$.
Thus we get a pointed $\CO_U$-module $\CZ=(\CO|_U)_0$.
We extend this module by zero outside the open set $U$ to obtain a pointed $\CO_X$-module which we likewise denote by $\CZ$.
We define a morphism $\al\colon\CZ\to \CF$ as follows.
For $V\subset U$ open, $\al_V\colon\CZ(V)\to\CF(V)$ is defined as $\al_v(a)=as_V$ for $a\in\CO_V$ and $\al_V(\omega_V)=0$.
This defines a morphism $\al$ in $\Mod_0(X)$.
Using $t$ instead of $s$ we define $\beta\colon\CZ\to\CF$ in the same manner.
Then $f\al=f\beta$ and since $f$ is a monomorphism, $\al=\beta$, hence $s_U=t_U$ and so $s_x=t_x$.
The other direction is trivial.

(b) This assertion is shown in the same way as for sheaves of abelian groups.
\qed

\subsection{Injectives and ascent}
\begin{proposition}
The category $\Mod_0(X)$ has enough injectives.
In particular, the class $\CI$ of products of skyscraper sheaves with injective stalks is an injective class.
\end{proposition}

\prf
Let $\CF$ be a pointed $\CO_X$-module.
For each point $x\in X$ the stalk $\CF_x$ is a pointed $\CO_{X,x}$-module.
Therefore there is an injection $\CF_x\hookrightarrow I_x$ into an injective $\CO_{X,x}$-module. 
Let $i_x$ denote the injection of $x$ in $X$ and consider the sheaf $\CI=\prod_{x\in X} i_{x,*}I_x$, which is a product of skyscraper sheaves with injective stalks.
For any $\CO_X$-module $\CG$ we have 
$$
\Hom_{\CO_X}(\CG,\CI)\cong \prod_{x}\Hom_{\CO_{X}}(\CG,i_{x,*}I_x)
$$ 
and for every $x\in X$ also 
$$
\Hom_{\CO_{X}}(\CG,i_{x,*}I_x)\cong \Hom_{\CO_{X,x}}(\CG_x,I_x).
$$
So there is a monomorphism $\CF\hookrightarrow\CI$ obtained from the maps $\CF_x\hookrightarrow I_x$.
It follows that $\CI$ is injective and hence the claim.
\qed

Let $\CC$ be the category of all sheaves of $\Q$-vector spaces on $X$. 
Consider the functor
$$
A\colon \Mod_0(X)\to \CC
$$
maps a sheaf $\CF$ to the sheafification of the presheaf 
$$
U\ \mapsto\  \Q[\CF(U)]/\Q x_0(U),
$$
where $x_0(U)$ is the special point of $\CF(U)$.

\begin{proposition}\label{2.8}
The functor $A$ is an $\CI$-injective ascent functor.
\end{proposition}

\prf
Since the ascent functor $A$ maps products of skyscraper sheaves to products of skyscraper sheaves the claim follows.
\qed

\subsection{Sheaf cohomology}
\begin{definition}
Let $X$ be a monoided space.
A sheaf $\CF$ is called \e{flabby} if for any two open sets $U\subset V$ the restriction map $\CF(V)\to \CF(U)$ is surjective.
\end{definition}

\begin{lemma}\label{lem2.5.2}
Every injective sheaf is flabby.
\end{lemma}

\prf
For any open set $U\subset X$ let $\CO_U$ denote the sheaf $j_!(\CO_X|_U)$, which is the restriction of $\CO_X$ to $U$, extended by zero outside $U$.
Now let $I$ be an injective $\CO_X$-module and let $U\subset V$ be open sets.
We have an inclusion $\CO_U\hookrightarrow \CO_V$ and since $I$ is injective we get a surjection $\Hom(\CO_V,I)\to\Hom(\CO_U,I)$.
But $\Hom(\CO_V,I)\cong I(V)$ and $\Hom(\CO_U,I)\cong I(U)$, so $I$ is flabby.
\qed

We consider the global sections functor $\Ga(X,\cdot)$ from $\Mod_0(X)$ to $\Mod_0(A)$.

\begin{lemma}
The global sections functor $\Ga(X,\cdot)$ is left strong-exact, admits an ascent datum, and sends injective maps to injective maps.
\end{lemma}

\prf
The left strong-exactness of the global sections functor is a standard verification.
Let $A$ denote the ascent functor above and let $A'$ be the ascent functor on $\Mod_0(\CO_X(X))$ given in Proposition \ref{2.2.2}.
It is easy to see that $A'\Ga=\Ga A$ holds on the full subcategory of $\CB=\Mod_0(X)$ consisting of flabby sheaves.
\qed

We define the cohomology of a sheaf $\CF\in\Mod_0(X)$ by
$$
H^i(X,\CF)\df R^i\Ga(X,\CF),\qquad i=0,1,\dots
$$

\begin{theorem}
The class of flabby sheaves is an acyclic class for the section functor $\Ga$.
In particular, we have
\begin{enumerate}[\rm (a)]
\item Every injective sheaf is flabby.
\item Let $0\to \CF\stackrel{f}{\to}\CH\stackrel{h}{\to}\CG\to 0$ be a strong exact sequence in $\Mod_0(X)$. If $\CF$ is flabby, then for every open set $U\subset X$ the sequence
$$
0\to \CF(U)\stackrel{f_U}{\to}\CH(U)\stackrel{h_U}{\to}\CG(U)\to 0
$$
is exact.
\item If in the situation of (a), the sheaves $\CF$ and $\CH$ are flabby, then so is $\CG$.
\item If $\CF$ is a flabby sheaf in $\Mod_0(X)$, then $H^i(X,\CF)=0$ for $i>0$.
\end{enumerate}
By Theorem \ref{thm1.6.11} we conlude cohomology can be computed using flabby resolutions.
\end{theorem}

\prf
(a) is Lemma \ref{lem2.5.2}.
(b) and (c) follow, after applying the ascent functor,  from the corresponding result for sheaves of abelian groups \cite{Hartshorne}.
Part (d) also follows after applying the ascent functor, since the ascent functor given maps flabby sheaves to flabby sheaves.
\qed

\begin{lemma}\label{3.9}
Let $\For$ be the forgetful functor from the category $\Mod_0(A)$ to $\Set_0\cong\Mod_0(1)$.
Then the isomorphy class of $\For(H^i(X,\CF))$ in $\Set_0$ does not depend on the choice of the sheaf $\CO_X$.
\end{lemma}

\prf
Let $\Set_0(X)$ denote the category of pointed sheaves over $X$ for the trivial structure sheaf $\CO_X=const$.
To compute the cohomology, use flabby resolutions in $\Mod_0(X)$.
They will remain flabby in $\Set_0(X)$.
\qed

\subsection{Noetherian Spaces}
\begin{definition}
We say that a monoid $A$ is \e{noetherian} if every chain of ideals $I_1\subset I_2\subset I_3\subset\dots$ is eventually stationary, i.e., there exists an index $j_0$ such that $I_j=I_{j_0}$ for every $j\ge j_0$.
A topological space $X$ is called \e{noetherian} if every sequence of closed subsets $Y_1\supset Y_2\supset Y_3\supset\dots$ is eventually stationary.
The \e{dimension} of a topological space is the supremum of the lengths of strictly descending chains of closed subsets.
A noetherian topological space is not necessarily of finite dimension.

If $X=\Spec\F_A$, then $X$ is noetherian if and only if $A$ is.
A monoid scheme $X$ is called \e{noetherian} if $X$ can be covered by finitely many affine schemes $\Spec(A_i)$ where each monoid $A_i$ is noetherian.
A noetherian scheme is noetherian and finite dimensional as topological space.
\end{definition}

Let $(\CF_\al)$ be a direct system of pointed sheaves.
By $\ds\lim_\to \CF_\al$ we denote the sheafification of the pre-sheaf $\ds U\mapsto \lim_\to\CF_\al(U)$.
Let $X$ be a monoided space.

\begin{lemma}
Let $(\CF_\al)_{\al\in I}$ be a direct system of flabby sheaves and assume that $X$ is noetherian.
Then $\ds \lim_\to\CF_\al$ is flabby.
\end{lemma}

\prf
As in the group valued case one proves that if $X$ is noetherian, then the pre-sheaf $\ds U\mapsto \lim_\to \CF_\al(U)$ already is a sheaf.
For every $\al\in I$ and every inclusion $V\subset U$ of open sets the restriction $\CF_\al(U)\to\CF_\al(V)$ is surjective.
This implies that $\ds\lim_\to\CF_\al(U)\to\lim_\to\CF_\al(V)$ also is surjective.
Since $X$ is noetherian we have $\ds\lim_\to\CF_\al(U_=(\lim_\to\CF_\al)(U)$, so $\ds\lim_\to\CF_\al$ is flabby.
\qed

Let $Y$ be a closed subset of $X$ and $\CF$ a pointed sheaf on $Y$.
Let $j_*\CF$ be the extension by zero outside $Y$.
Then one has $H^i(Y,\CF)=H^i(X,j_*\CF)$ as a flabby resolution $\CJ^\bullet$ of $\CF$ gives a flabby resolution $j_*\CJ^\bullet$ of $j_*\CF$.

\begin{theorem}
Let $X$ be noetherian of dimension $n$, and let $\CF$ be a pointed sheaf which is generated by finitely many sections.
Then for every $i>n$ we have $H^i(X,\CF)=0$.
\end{theorem}

\prf
By Lemma \ref{3.9} we may assume that $\CO_X$ is the trivial sheaf of monoids.
For a closed subset $Y$ of $X$ and a pointed sheaf $\CF$ on $X$ we write $\CF_Y$ for $j_*(\CF|_Y)$.
If $U\subset X$ is open, we write $\CF_U= i_!(\CF|_U)$.
Then, if $U-X\setminus Y$, we have an exact sequence
$$
0\to \CF_U\hookrightarrow \CF\to \CF_Y\to 0,
$$
as one easily checks.

We next reduce the proof to the case when $X$ is irreducible.
For assume $X$ is reducible, then $X=Y\cup Z$ with closed sets $Y, Z$ both different from $X$.
Let $U=X\setminus Y$ and consider the exact sequence
$$
0\to\CF_U\hookrightarrow\CF\to\CF_Y\to 0.
$$
By the long exact sequence of cohomology it suffices to show $H^i(X,\CF_U)=0$ and $H^i(X,\CF_Y)=0$.
Now $\CF_U$ can be viewed as a sheaf on $Z$ and so the proof if reduced to the components $Y$ and $Z$.
By induction on the number of components we can now assume that $X$ is irreducible.

We prove the Theorem by induction on $n=\dim X$.
If $n=0$ then $X$ has only two open sets, itself and the empty set.
Then $\Ga(X,\cdot)$ is exact, so the claim follows.
Now for the induction step let $X$ be irreducible of dimension $N>0$ and let $\CF$ be a pointed sheaf on $X$.
By an induction argument it suffices to assume that $\CF$ is generated by a single section in $\CF(U)$, say, for an open set $U$.
Let $\CZ$ be the constant sheaf with fibre $\Z/2\Z$.
Then $\CF$, being generated by a single section, is a quotient of $\CZ_U$.
So we have an exact sequence,
$$
0\to \CR\hookrightarrow \CZ_U\to\CF\to 0.
$$
By the long exact cohomology sequence it suffices to show the vanishing of the cohomology of $\CR$ and $\CZ_U$.
If $\CR\ne 0$, then there exists an open set $V\subset U$ such that $\CR_V\cong \CZ_V$.
So we have an exact sequence
$$
0\to \CZ_V\hookrightarrow\CR\to\CR /\CZ_V\to 0.
$$
The sheaf $\CR/\CZ_V$ is supported in $\ol{U\setminus V}$ which has dimension $<n$ since $X$ is irreducible.
So it follows that $H^i(X,\CR /\CZ_V)=0$ for $i>n$ by induction hypothesis.
It remains to show vanishing of cohomology for $\CZ_V$.
We show that for every open $U\subset X$ we have $H^i(X,\CZ_U)=0$ for $i>n$.
Let $Y=X\setminus U$.
We have an exact sequence
$$
0\to \CZ_U\hookrightarrow \CZ\to \CZ_Y\to 0.
$$
Since $X$ is irreducible, we have $\dim Y<n$.
So by induction hypothesis we have $H^i(X,\CZ_Y)=0$ for $i\ge n$.
On the other hand, $\CZ$ is flabby as it is a constant sheaf on an irreducible space.
Hence $H^i(X,\CZ)=0$ for $i>0$.
So the long exact cohomology sequence gives the claim.
\qed

\subsection{Base change}
Now assume that $X$ is a monoid scheme.
Let $X_\Z= X\otimes\Z$ be the base change to $\Z$.
Instead of $\Z$ one could take any other ring here.
Let $\CF$ be a pointed sheaf over $X$.
For a pointed module $(M,m_0)$ over a monoid $A$ write $M_\Z$ for the $\Z[A]$-module $\Z[M]/\Z m_0$.
Every open set $U$ in $X$ defines an open set $U_\Z$ in $X_\Z$ as follows.
If $X=\Spec_{\F_1}A$ is affine, then $U$ defines an ideal $\a$ of $A$.
Then $\Z[\a]$ is an ideal of $\Z[A]$ which defines an open set $U_\Z$ of $X_\Z=\Spec\Z[A]$.
For non-affine $X$ define $U_\Z$ locally and take the union.
We define the sheaf $\CF_\Z$ to be the sheafification of the pre-sheaf
$$
U\ \mapsto\ \lim_{\stackrel{\longrightarrow}{V_\Z\supset U}}\CF(V)_\Z.
$$
here the inductive limit is taken over all open sets in $X_\Z$ which contain $U$ and are of the form $V_\Z$ for some $V$ open in $X$.

If $\CF$ is a skyscraper sheaf in $x\in X$, then the closed set $\bar x=\ol{\{ x\} }$ is given by an ideal sheaf which base changes to an ideal sheaf of $X_\Z$ which defines a closed subset $\bar x_\Z$ of $X_\Z$.
It turns out that $\CF_\Z$ is a constant sheaf on $\bar x_\Z$ extended by zero outside $\bar x_\Z$.
In particular, $\CF_\Z$ is flabby.

The functor $\CF\mapsto\CF_\Z$ is an ascent functor from $\Mod_0(X)$ to $\Mod( X_\Z)$ which maps sheaves in the injective class $\CI$ to flabby sheaves, hence $\CI$-resolutions are mapped to flabby resolutions.

\begin{theorem}
There is a natural injection,
$$
H^p(X,\CF)_\Z\ \hookrightarrow\ H^p(X_\Z,\CF_\Z).
$$
\end{theorem}

\prf
Let $0\to \CF\to I^0\to I^1$ be an injective resolution, where $I^p$ is a product of skyscraper sheaves.
Lemma \ref{ascentProp} gives an injection $H^p(\Ga I)_\Z\hookrightarrow H^p((\Ga I)_\Z)$.
As $I$ consist of product of skyscraper sheaves, the complex $(\Ga I)_\Z$ is isomorphic with $\Ga(I_\Z)$.
As $I_\Z$ is a flabby resolution of $\CF_\Z$, we get
$$
H^p(X,\CF)_\Z=H^p(\Ga I)_\Z\hookrightarrow\ H^p(\Ga(I_\Z))=H^p(X_\Z,\CF_\Z).
$$
\qed

\begin{corollary}
If $X=\Spec_{\F_1}(A)$ is affine and $M$ is a pointed $A$-module, then $H^p(X,\tilde M)=0$ for $p>0$.
\end{corollary}

\prf
Since $(\tilde M)_\Z\cong \widetilde{M_\Z}$, the claim follows from the corresponding claim for schemes.
\qed

\today

{\small Mathematisches Institut\\
Auf der Morgenstelle 10\\
72076 T\"ubingen\\
Germany\\
\tt deitmar@uni-tuebingen.de}

\end{document}